\documentclass[12pt, reqno]{amsart}
\usepackage{amsmath, amstext, amsbsy, amssymb, amscd}

\newtheorem{theorem}{Theorem}[section]
\newtheorem{lemma}[theorem]{Lemma}

\newtheorem{corollary}[theorem]{Corollary}
\theoremstyle{definition}
\newtheorem{definition}[theorem]{Definition}
\newtheorem{example}[theorem]{Example}

\theoremstyle{remark}
\newtheorem{remark}[theorem]{Remark}

\numberwithin{equation}{section}

{\vskip-\lastskip\medskip
  \noindent
  {\em #1.}\enspace
  }%
{\qed\par\medskip
  }

\begin{document}

\title
[ Affine  Algebraic     Varieties]
{  Affine  Algebraic     Varieties  }

\author[Jing  Zhang]{Jing  Zhang}
\address{ Department of Mathematics and Statistics, 
University at  Albany, SUNY, Albany, NY 12222,  USA}
\email{jzhang@albany.edu}

\begin{abstract} In this paper, 
we  give     new criteria 
for  affineness of  a variety  
defined  over   $\Bbb{C}$. 
Our  main result is that
 an irreducible algebraic   variety  $Y$  
 (may be singular)
 of dimension $d$ ($d\geq 1$) defined over 
$\Bbb{C}$ 
is an  affine  variety
 if and only if $Y$  contains no complete curves, 
$H^i(Y, {\mathcal{O}}_Y)=0$  for all $i>0$ 
and     the boundary   $X-Y$  is  support of a big  divisor, 
 where $X$  is  a   projective variety 
containing
$Y$.
We  construct   three examples
to show   that a variety is not affine if it 
only satisfies  two conditions   among  these  three conditions.  
We also give examples to demonstrate  the difference  
between the  behavior  of the boundary divisor
$D$  and  the affineness of  $Y$.

 If
$Y$ is an affine variety, then 
the ring  $\Gamma (Y, {\mathcal{O}}_Y)$
is noetherian. 
However, to prove that $Y$  is an affine variety,
we do not start  from  this ring. 
We explain why we do not need to check 
the noetherian property 
of the ring  $\Gamma (Y, {\mathcal{O}}_Y)$ directly
but use 
 the techniques of sheaf and cohomology.

\end{abstract} 

\maketitle

\begin{flushleft}
\end{flushleft}
\date{}
\section{Introduction}

We work over complex number field 
$\Bbb{C}$. 

Affine varieties  are  important  in   algebraic  geometry. 
In    1957, 
 J.-P. Serre 
 introduced sheaf and cohomology techniques 
 to algebraic geometry  and 
  discovered  his well-known cohomology 
criterion (\cite{26}, \cite{8}, Chapter 2, Theorem 1.1): 
a variety  (or  a  noetherian  scheme) $Y$ is affine 
if and only  if  for all coherent 
sheaves $F$  on $Y$  and all  $i>0$, $H^i(Y, F)=0$.
The point of Serre's   criterion
is that instead of  to  look  at the noetherian  property of
the ring 
$\Gamma (Y, {\mathcal{O}}_Y)$, 
  to check the  affineness of  $Y$,
  we examine   the cohomology of 
  the coherent sheaves on $Y$. 
  In 1962, Narasimhan solved the Levi
  problem for complex spaces (\cite{21}, \cite{22}). 
  The corresponding analytic variety in
   complex geometry,   which holds the similar important  position
   as an affine  variety in algebraic geometry,   is  a Stein variety.
   We know that to verify  the Steinness of an
   analytic variety  $Y$, we do not look at the 
   ring  $\Gamma (Y, {\mathcal{O}}_Y)$, 
   but consider the holomorphic functions and
   check  whether  $Y$  is holomorphically separable and
   holomorphically convex  (\cite{6}, Page 143).  
     Inspired by Serre's  criterion and the analytic method of 
  Narasimhan,  in 1969, 
Goodman and Hartshorne 
proved that $Y$ is an affine variety if and only if 
$Y$ contains no complete curves and the dimension $h^1(Y, F)$
of the linear space 
$H^1(Y, F)$ is bounded for all coherent sheaves $F$ on $Y$ \cite{4}.

Let $X$ be the completion of $Y$.  In 1969, 
Goodman also proved that $Y$ is affine if and only if 
after suitable blowing up
the closed subvariety on the  boundary  $X-Y$, the new boundary $X'-Y$
is  support of an ample divisor, where $X'\rightarrow X$
is the  blowing up with center in $X-Y$
(\cite{3}; \cite{8}, 
 Chapter 2, Theorem 6.1).  
 For any quasi-projective
 variety  $Y$, we   may  assume that 
 the boundary $X-Y$  is the support of an effective
 divisor  $D$  with simple normal crossings
 by blowing up the closed subvariety in $X-Y$. 
   $Y$ is affine if $D$ is
ample. So  if we can show the ampleness of
$D$, $Y$ is affine. There are two important 
criteria for ampleness due to Nakai-Moishezon and 
Kleiman  
 (\cite{12}; \cite{15}, Chapter 1, Section  1.5). 
Another sufficient condition is that if 
$Y$ contains no complete curves and the linear system 
$|nD|$ is base point free, then $Y$ is  affine  
(\cite{8}, 
Chapter 2, Page 64).  Therefore we can apply base point free theorem
if we know the numerical condition of $D$ (\cite{24},  \cite{15}, 
Chapter 3, Page 75, Theorem 3.3). 
In 1988, 
Neeman proved that if $Y$
 can be embedded in an affine 
scheme Spec$A$, then $Y$ is affine if and only if 
$ H^i(Y, {\mathcal{O}}_Y )=0 $  for all $i>0$  \cite{23}.
The 
significance of 
 Neeman's Theorem   is that  it  is  not assumed that 
the ring $A$ is noetherian.

In higher dimension (at least, in our problem),
it is very hard to check  the ampleness of a big
(even big and nef)  divisor $D$
and  the base point freeness of the linear system 
$|nD|$. We  need to search for a different approach.

 Iitaka's $D$-dimension theory is widely 
used in classification  of algebraic varieties
(\cite{10}, \cite{11} \cite{19}, \cite{30}). Recall the definition 
of notation $\kappa(D, X)$. 

\begin{definition}
 If for all $m>0$, $H^0(X, {\mathcal{O}}_X(mD))=0$, then 
we define the $D$-dimension  (or  Iitaka dimension)
$\kappa(D, X)=-\infty$. Otherwise, 
$$\kappa(D, X)=tr.deg_{\Bbb{C}}\oplus_{m\geq 0}
H^0(X, {\mathcal{O}}_X(mD))-1.$$
\end{definition}

In particular, 
the Kodaira dimension $\kappa(X)$ of   $X$ 
is defined to be   $\kappa(K_X, X)$,
where  $K_X$  is the canonical divisor of  $X$. 

A  quasi-affine variety  is  a Zariski 
open subset of an affine variety.
Throughout  this paper, we assume that 
$Y$ is   an  irreducible  algebraic  variety  of  dimension $d$  
defined  over  $\Bbb{C}$
such that  
it is  the complement of an effective  divisor  $D$  in a 
projective  variety  $X$.  We may assume that $D$  has simple normal 
crossings.

\begin{theorem} An  irreducible  algebraic surface $Y$
is an affine surface if and only if 
$Y$  contains no complete curves, the boundary 
$X-Y$ is connected and 
 $\kappa(D, X)=2$. 

\end{theorem}

Theorem 1.2 is not true for higher dimensional varieties 
(see the example in Section 3).

\begin{theorem} An  irreducible   algebraic  variety  $Y$  of dimension
 $d$ ($d\geq 1$)
is a  quasi-affine    variety  if and only if $Y$  contains no complete curves
 and  
$\kappa(D, X)=d$. 
\end{theorem}

\begin{theorem} An  irreducible   algebraic  variety  $Y$  of dimension
 $d$ ($d\geq 1$)
is an  affine  variety  if and only if $Y$  contains no complete curves, 
$H^i(Y, {\mathcal{O}}_Y)=0$  for all $i>0$ and  
$\kappa(D, X)=d$. 
\end{theorem}

If we have a surjective morphism from a variety $Y$ to
an affine variety such that every fibre is affine,
then $Y$  may not be affine  (see  Example 3.11 in Section 3). 
But we have the following theorem.

\begin{theorem} Let $f: Y\rightarrow W$  be
a surjective morphism  from an 
irreducible   quasi-projective
variety  $Y$  to an  irreducible    affine  variety 
$W$. If every fiber is an affine subvariety 
of  $Y$  and    
$H^i(Y, {\mathcal{O}}_Y)=0$  for all $i>0$,
then  $Y$  is  an affine variety. 
\end{theorem}

\begin{corollary}
 Let $f: Y\rightarrow W$  be
a surjective morphism  from an 
irreducible   quasi-projective
variety  $Y$  to an  irreducible    affine  variety 
$W$. If 
$Y$  has no complete curves, 
a general fiber is an affine subvariety 
of  $Y$  and    
$H^i(Y, {\mathcal{O}}_Y)=0$  for all $i>0$,
then  $Y$  is  an affine variety. 
\end{corollary}

In Section 3, we will give three examples to   
show   that 
a variety is not affine if it only satisfies two 
conditions in  Theorem  1.4.

 If  $Y$  is an affine
variety, then
the ring  $\Gamma(Y, {\mathcal{O}}_Y)$  is  noetherian. However, 
in our proof  of Theorem 1.4, we do not directly check the 
noetherian  property of this ring,
which I do not know  whether  it is possible or not,
because our condition is rather geometric. 
 In  \cite{36}, we proved Mohan Kumar's affineness 
 conjecture for an algebraic manifold 
 and gave a partial answer to J.-P. Serre's  Steinness
 question \cite{27}. 
 After we carefully   examined the conditions 
 and where we used them,
 we found that we can get general theorems
 for  singular varieties  by changing 
 the assumption and modifying the proof.  The vanishing 
 Hodge cohomology  $H^i(Y, \Omega^j_Y)=0$
 for all $i>0$  and $j\geq 0$  
 are replaced by two conditions:
 $Y$  has no complete curves  and 
 $H^i(Y, {\mathcal{O}}_Y)=0$  for all 
 $i>0$.  The advantage of these two conditions
 is that  they are well-defined 
 for singular varieties. 
 This makes the generalization  possible.

The question  of  a quasi-projective variety  $Y$ to be affine  is
different from the 
behavior  of the boundary divisor $D$, in particular, 
the numerical condition of $D$ like nefness  and  
finitely generated  property of  the graded  ring
$$\oplus_{n= 0}^\infty H^0(X, {\mathcal{O}}_X(nD)).$$ 
The reason is that 
$$\Gamma(Y, {\mathcal{O}}_Y )
\neq \oplus_{n= 0}^\infty H^0(X, {\mathcal{O}}_X(nD)).$$

 We will give two examples to 
demonstrate  this difference  in Section 3. One example (due to Zariski)
 gives a   surface   $Y=X-D$,   which  is  affine   but the    
the corresponding  graded  ring
$\oplus_{n=0}^{\infty} H^0(X, {\mathcal{O}}_X(nD))$  is not  finitely generated
for  an  effective    divisor  $D$. The other example (\cite{ 16}; \cite{8}, Page 232) is
a surface $Y=X-D$  such that  
$$H^0(Y, {\mathcal{O}}_Y)=H^0(X, {\mathcal{O}}_X(nD))=\Bbb{C}.$$

A necessary condition for the affineness of  $Y$   with dimension  $d$  is that 
$Y$  must have  plenty of  nonconstant regular functions.
More  precisely, 
 $Y$  must have  $d=$dim$Y$  algebraically independent nonconstant
regular functions.  This  means that   the corresponding effective  
boundary  divisor
$D$  must be big,   i.e.,    
$$ h^0(X, {\mathcal{O}}_X(nD))\geq  a n^d
$$ 
for some positive  number  $a$  and $n\gg 0$. 
So  the above surface  $Y$  is not affine  but 
$\oplus_{n= 0}^\infty H^0(X, {\mathcal{O}}_X(nD))$
is  finitely  generated.

We will prove the results in Section 2 and give
examples in Section 3.

\section{Proof of the Theorems}

Recall our notation:  $Y$ is an open subset of 
a projective variety $X$   with  dimension $d\geq 1$ 
and $D$ is the effective boundary divisor
 with support $X-Y$.  We may assume that the boundary 
divisor $D$ has simple normal crossings by further blowing up suitable
closed subvariety  of $X-Y$. 

{\bf Proof of Theorem 1.2.} 
The proof is the same as the proof of  Theorem 1.1 \cite{35}.
The idea is to show  that 
the boundary 
$X-Y$ is support of an ample divisor. 
\begin{flushright}
 Q.E.D. 
\end{flushright}

{\bf Proof of Theorem 1.3.} 
One direction is trivial. 
If $Y$ is an  irreducible  quasi-affine  variety,
 then
it has no
complete curves since  it is a closed   subset of 
affine space  ${\Bbb{C}}^N$. 
Let  $U$  be  an affine  variety 
such that  $Y\subset  U \subset  X$, 
then the boundary $X-U$  is of pure 
codimension 1   (\cite{8}, Chapter II, Proposition 3.1).
By Goodman's theorem (\cite{8}, Chapter II, Theorem 6.1),
since  $U $  is affine 
we may assume that  $X-U$  is support 
of an ample divisor 
$A$  on  $X$. Therefore 
$\kappa(A, X)=d$ (\cite{15}, Proposition 
2.61). 
The support of   $A$
is contained in the support of  $D\geq 0$,
so   $\kappa(D, X)=\kappa(A, X)=d$.

Suppose now that 
$Y$  satisfies  the condition in
Theorem 1.3. We will prove that it is quasi-affine. 
We may assume that both 
$X$ and  $Y$  are normal. In fact, let $Y'$
be the normalization of  $Y$, then 
we have finite morphism 
from  $Y'$  to  $Y$. Thus 
$Y'$  
satisfies  two conditions in  Theorem  1.3
and  $Y$  is  affine if  and  only if
$Y'$  is  affine.

\begin{lemma} Theorem  1.3  holds for curves   and  surfaces. 

\end{lemma}

$Proof.$ If a curve  is not complete then it is affine  (\cite{8}, Chapter II, 
Proposition 
4.1).  Surface case is  true  by Theorem 1.1. 
\begin{flushright}
 Q.E.D. 
\end{flushright}

We may
assume that Theorem 1.3  holds for 
$(d-1)$-dimensional   algebraic  varieties.
Suppose  dim$Y=d$. 

 \begin{lemma} Under the condition of  Theorem 1.3,
 if  $Y$  is smooth, then  for every  point 
 $y\in  Y$, there is  a smooth prime principle 
 divisor $Z$  passing through 
 $y$. 
\end{lemma}

$Proof$. Since $X$ is projective  and smooth, there is a hypersurface 
$H_1$  defined by a homogeneous  polynomial $h_1$ of degree  at  least 
2  passing through $y_1$ and  $Z=H_1\cap X$  is a prime   principal 
 divisor on
$X$.

Let $H_2$   be  a distinct hypersurface  defined by 
$h_2$ with the same degree
such that   $h_2(y_1)\neq 0$ and 
$H_2\cap X$  is a prime principal divisor on
$X$. Let  $h=\frac{h_1}{h_2}$, then  $h$
is a rational function on $Y$ and regular on an open subset of $Y$ containing 
$y_1$.  By the following claim,  there are two regular functions $f$  and $g$ on
$Y$  such that   $h=\frac{h_1}{h_2}=\frac{f}{g}$. 
So $h_1=0$  if and only if  $f=0$. Therefore the subvariety   
$Z$ in $Y$ is defined by $f$. 

If  $h^0(X, {\mathcal{O}}_X(mD))>  0$
for  some  $m\in {\Bbb{Z}}$  and $X$  is normal, 
choose a basis $\{f_0, f_1, \cdot \cdot\cdot, f_n\}$
     of the linear space 
     $H^0(X, {\mathcal{O}}_X(mD))$, it defines a rational 
     map 
     $\Phi _{|mD|}$
     from $X$ to the projective space 
     ${\Bbb{P}}^N$ by sending a point $x$ on $X$ to
     $(f_0(x), f_1(x), \cdot \cdot\cdot, f_N(x))$ in ${\Bbb{P}}^N$. 
By definition of $D$-dimension \cite{30}, Definition 5.1,
$$ \kappa (D, X)= \max_m\{\dim (\Phi _{|mD|}(X))\}. 
      $$

Let ${\Bbb{C}}(X)$ be the function field of $X$. Let 
$$  R(X, D)=\oplus_{\gamma= 0}^\infty H^0(X, {\mathcal{O}}_X(\gamma D)) 
$$
be the graded $\Bbb{C}$-domain  and $R^*\subset R$
the  multiplicative  subset  of all nonzero
homogeneous  elements. Then the quotient ring $R^{*-1}R$  
is a graded  $\Bbb{C}$-domain  and its degree 0 part 
$(R^{*-1}R)_0$  is a  field  denoted  by   $Q((X, D))$, i.e., 
$$ Q((X, D))=(R^{*-1}R)_0.
$$  

Let $X$ be  normal proper  over 
$\Bbb{C}$, then we have \cite{19}.

(1) If there is an $m_0>0$  such that  for all $m>m_0$, 
$h^0(X, {\mathcal{O}}_X(mD))>0$, then 
$$ {\Bbb{C}}(\Phi _{|mD|}(X)) = Q((X, D)).
$$  

(2) If  $\kappa(D, X)=$\mbox{dim}X, then  
$\Phi _{|mD|}$ is birational for all $m\gg 0$. 
In particular, ${\Bbb{C}}(X)=Q((X, D))$. 

So if $D$  is a big divisor, then 
any rational function on  $X$
can be written as a quotient of 
two elements in 
$H^0(X, {\mathcal{O}}_X(mD))$  for sufficiently 
large  $m$. These two elements are regular 
on  $Y$.

In the proof, we used the generalized version of Bertini's Theorem: 
if  $X$  is a  projective  algebraic 
manifold of dimension at least 2, then  for any point  $x_0$  on  $X$,
there  is  an  irreducible     smooth  hypersurface 
$H$  of  degree at least 2 passing through  $x_0$   such that  $H$
intersects   $X$  with  an irreducible  smooth  codimension 1  subvariety  
of  $X$  \cite{36}. We give a proof here for 
completeness. 

 Let  $X$  be  a closed subset of  ${\Bbb{P}}^n$,
$n\geq  3$.
We may assume  that  $X$  is 
not contained in any hyperplane   and  after coordinate transformation,
the homogeneous coordinate of  $x_0$
is  $(1,0,...,0)$.

Let $H$  be a    hypersurface  defined by a homogeneous 
polynomial  $h$   of   degree   2
passing  through  $x_0$, then 
$$h=\sum_{j=1}^na_{0j}x_0x_j
+\sum_{i=1}^n\sum_{j\geq  i}a_{ij}x_ix_j.
$$  
$H$  is nonsingular  at 
$x_0$  if at least one $a_{0j}\neq  0$. 
Let  $V$  be the linear space of these hypersurfaces, then
the dimension of  $V$  is  
$$\mbox{dim}_{\Bbb{C}}V=
\frac{(n+2)(n+1)}{2}-1=\frac{n^2+3n}{2}.$$

By  Euler's formula, the hypersurface is singular at a point
 $x=(x_0, x_1,..., x_n)$ if and only if 
$$\frac{\partial h}{\partial x_0} =\frac{\partial h}{\partial x_1} =...=
\frac{\partial h}{\partial x_n}=0. 
$$
It is    a  system of  linear  equations

     $$\quad\quad\quad \quad  a_{01}x_1 + a_{02}x_2 + \cdots    + a_{0n}x_n  =0$$
    $$ a_{01}x_0 + 2a_{11}x_1 + a_{12}x_2 + \cdots    +  a_{1n}x_n  =0$$
    $$\cdot\cdot\cdot\cdot\cdot\cdot\cdot\cdot\cdot\cdot\cdot\cdot\cdot
\cdot\cdot\cdot\cdot\cdot\cdot\cdot\cdot\cdot
\cdot\cdot\cdot\cdot\cdot\cdot\cdot\cdot\cdot\cdot\cdot\cdot$$ 
    $$a_{0n}x_0       + a_{1n}x_1 +a_{2n}x_2 +\cdots    +  2a_{nn}x_n  =0$$

The above system has a solution in  ${\Bbb{P}}^n$
if and only if  the determinant   of  
the following 
 symmetric  matrix  $A$ is zero,

\[\left( \begin{array}{llll}
 0  & a_{01}  & a_{02} \cdots & a_{0n}\\
a_{01} & 2a_{11} & a_{12} \cdots & a_{1n}\\
\multicolumn{4}{c}\dotfill\\
a_{0n} & a_{1n}  & a_{2n} \cdots  & 2a_{nn}
\end{array}   \right). \]

Considering  $(a_{01}, a_{02},..., a_{(n-1)n})$
as a point in the projective space  ${\Bbb{P}}^{\frac{n^2+3n}{2}-1}$,
the system has solution only on the hypersurface  defined  by
det$A=0$.  So  the degree 2  hypersurface $H$ in $V$  is nonsingular  on an 
open subset of  ${\Bbb{P}}^{\frac{n^2+3n}{2}-1}$, i.e., a general 
member  $H$   of   $V$  is smooth.

Let  $x$  be a   closed  point of  $X$  and  define  $S_x$  to  be the set of
smooth hypersurfaces  $H$  (defined by  $h$) of degree
2
such that      $x$  is  a  singular 
point of  $H\cap  X$. $H$  is ample and $X$  is not contained in $H$  by our assumption so  $X\cap H$ 
is not empty.  Fix  a smooth  irreducible  hypersurface  $H_0$
of degree 2 such  that  $x$  is not a point of  $H_0$. Let $h_0$
be the  defining   homogeneous polynomial of $H_0$, then
$h/h_0$  gives  a  regular function 
 on  ${\Bbb{P}}^n-H_0$.  
When restricted on  $X$, it is  a regular 
function on  $X-X\cap  H_0$.

  Define   a map $\xi_x$   from  the  linear space   $V$   to 
${\mathcal{O}}_{x, X}/{\mathcal{M}}_x^2$    as follows:
for every  element $h$  in  $V$   (a homogeneous  polynomial 
of degree 2 such that the  corresponding  hypersurface  $H$ is smooth and 
passes through  $x_0$),  the image  $\xi_x(h)$   is the image  
of  $h/h_0$   in the local  ring  
${\mathcal{O}}_{x, X}$  
 modulo   ${\mathcal{M}}_x^2$.   It is easy to see that
$x$  is a point of   $H\cap   X$  if and only if  the 
image $\xi_x(h)$ of the defining polynomial  $h$ of  $H$ 
is  contained in   ${\mathcal{M}}_x$.  And
$x$ is singular   on  $H\cap   X$  if and only if  the 
image $\xi_x(h)$   is  contained in   ${\mathcal{M}}_x^2$
because  the local ring  ${\mathcal{O}}_x/\xi_x(h)$
will not be regular. 
So there is the following one-to-one  correspondence
$$H\in  S_x \Longleftrightarrow  h\in   ker  \xi_x.  
$$      
Since   $x$ is a closed point and  the  ground field is  
$\Bbb{C}$,  the maximal  ideal  
${\mathcal{M}}_x$    is   generated   by  linear  forms  in  the coordinates.
The map  $\xi_x$  is surjective if $x\neq  x_0$.
The map  $\xi_{x_0}$  is not surjective to 
${\mathcal{O}}_{x, X}/{\mathcal{M}}_x^2$ 
but    
surjective  to  
${\mathcal{M}}_{x_0}/{\mathcal{M}}^2_{x_0}$.

Let  $d$  be the dimension of   $X$,  then  the linear space  
${\mathcal{O}}_{x, X}/{\mathcal{M}}_x^2$  has dimension  $d+1$  
over  $\Bbb{C}$. Considering  the  map  
$$\xi_x: V\longrightarrow    {\mathcal{O}}_{x, X}/{\mathcal{M}}_x^2,
$$
if  $x\neq  x_0$, 
the kernel  has dimension 
$$  \mbox{dim}_{\Bbb{C}}ker \xi_x=\frac{n(n+3)}{2}-d-1. 
$$
Therefore  the linear space  $S_x$  is  a  linear  system  of
hypersurfaces    with
dimension   $\frac{n(n+3)}{2}-d-2$  if 
$x\neq  x_0$.  If  $x=x_0$, then  the  projective 
dimension
of  $S_{x_0}$  is   
 $\frac{n(n+3)}{2}-d-1$.

Considering  the  linear system  $V$  as a  projective space, 
then  $X\times  V$ is a projective  variety.  The subset  
$S\subset  X\times  V $   consists    of all   pairs   $<x, H>$
such that  $x\in  X$  is  a  closed  point  and  
$H\in   S_x$.

        $S$  is the set of  closed  points   of  a closed  subset   of
$  X\times  V $  and  we   give    a reduced  induced  scheme structure
to  $S$.  The  first  projection   
$p_1:  S\rightarrow   X$  is surjective. 
If $x\neq  x_0$,  
 the  fiber 
 $p_1^{-1}(x)$   
is a projective space   with  dimension   
$\frac{n(n+3)}{2}-d-2.$ 
The special fiber  $p_1^{-1}(x_0)$
is a projective space with
 dimension $\frac{n(n+3)}{2}-d-1.$ 
   Hence  $S$     has    
dimension   
$$ [\frac{n(n+3)}{2}-d-2]+d=\frac{n(n+3)}{2}-2.
$$

Let  $S=\cup_{i=0}^mS_i $  be
 an  irreducible  decomposition.
 Then  every  $p_1(S_i)$  
 is closed  and 
 there is an $i$, such that 
 $p_1(S_i)=X$.  For every  $S_i$   with
   $p_1(S_i)=X$,  there is an open subset 
   $U_i\subset S_i$ such that for every $x\in U_i$,
   the fiber $p_1^{-1}(x)$ has  constant dimension 
   $n_i$.  Let $x\in \cap U_i$, since  
   the fiber  $p_1^{-1}(x)$  is irreducible, it is contained in some $S_i$. Suppose $p_1^{-1}(x)\in S_1$.
   Let $f_1$ be the restriction of  $p_1$ on 
   $S_1$, i.e., $p_1|_{S_1}=f_1$, then  
   $p_1^{-1}(x)\subset  f_1^{-1}(x)$  since  $p_1^{-1}(x)$
   is  irreducible. 
   The opposite inclusion  is  obvious, so
   $p_1^{-1}(x)=f_1^{-1}(x)$ for $x\in \cap U_i$
   and  $n_1=\frac{n(n+3)}{2}-d-2$. 
   
   Since  $f_1$  is surjective and $S_1$  is one
   irreducible  component  of  $S$, for every 
   $x\in X$, the fiber  $f_1^{-1}(x)$  
   is  not empty and contained in 
    $p_1^{-1}(x)$.  But  the dimension of 
   $f_1^{-1}(x)$  is  at least 
   $\frac{n(n+3)}{2}-d-2$, so 
   for every 
   $x\in X$, $p_1^{-1}(x)=f_1^{-1}(x)$. 
   Thus  $S_1=S$  and  $S$  is irreducible.

Looking   at  the   second  projection    (a proper  morphism )
$p_2: S\rightarrow    V$.  The dimension of the image
$$  \mbox{dim}  p_2(S)\leq  \mbox{dim}  S=\frac{n(n+3)}{2}-2.
$$
Since  $S$  is closed   in   $  X\times  V $
and the dimension of $V$ (as a projective space)
is  $\frac{n(n+3)}{2}-1$,
$V-p_2(S)$  is an open   subset   of  
$V$. This implies that a general member $H$ of  $V$
intersects  $X$  
with a smooth  variety $X\cap  H$.

Next we will prove that  $X\cap  H$
is   irreducible.

From the short exact  sequence
$$ 0\longrightarrow 
 {\mathcal{O}}_{{\Bbb{P}}^n}(-H)
\longrightarrow 
 {\mathcal{O}}_{{\Bbb{P}}^n}
\longrightarrow 
 {\mathcal{O}}_{H}
\longrightarrow 
0,
$$
since  
$H^1({\Bbb{P}}^n, {\mathcal{O}}_{{\Bbb{P}}^n}(-H))=0$  
(\cite{7}, Page 225, Theorem 5.1),
we have   a  surjective map 
$$H^0({\Bbb{P}}^n,{\mathcal{O}}_{{\Bbb{P}}^n})=\Bbb{C}
\longrightarrow  H^0(H, {\mathcal{O}}_{H}).$$
So $H^0(H, {\mathcal{O}}_{H})=\Bbb{C}$  and 
 the hypersurface  $H$  is connected.

Since  $X$  is  closed in  ${{\Bbb{P}}^n}$, 
$H|_X$  is ample  on  $X$.   By 
Kodaira  Vanishing  Theorem,
$H^1(X, {\mathcal{O}}_{X}(-H))=0$  (\cite{15}, Page  62).
Applying the short exact sequence
$$ 0\longrightarrow 
 {\mathcal{O}}_X(-H)
\longrightarrow 
 {\mathcal{O}}_X
\longrightarrow 
 {\mathcal{O}}_{H\cap  X}
\longrightarrow 
0,
$$ 
we  get  
$$H^0(H\cap X,  {\mathcal{O}}_{H\cap  X})=H^0(X, {\mathcal{O}}_X)=\Bbb{C}.$$
Thus  the intersection  $H\cap  X$  is connected.  
Therefore  for a general hypersurface 
$H$  of degree 2,   $H\cap  X$  is   smooth  and  irreducible.

We have proved that  a general  smooth  hypersurface   of  degree  2  
 passing   through   
$x_0$   intersects   $X$  with  an  irreducible  
 smooth  subvariety   of codimension  1.  

\begin{flushright}
 Q.E.D. 
\end{flushright}

\begin{lemma} Under the condition of  Theorem 1.3, 
any prime principal divisor $Z=\{f=0, f\in H^0(Y, {\mathcal{O}}_Y)\}$ 
satisfies the same condition,
i.e., $Z$ contains no complete curves, 
$H^i(Z, {\mathcal{O}}_Z)=0$  for all $i>0$ and  
$\kappa(D|_{\bar{Z}}, \bar{Z})=d-1$, where  $\bar{Z}\subset X$ 
is a closed subvariety of $X$ containing $Z$.
\end{lemma}

$Proof.$ It is obvious that $Z$  contains no complete curves
since  $Z$  is   a  subvariety  of   $Y$.

From the short exact sequence
$$ 0\longrightarrow 
 {\mathcal{O}}_Y
\longrightarrow 
 {\mathcal{O}}_Y
\longrightarrow 
 {\mathcal{O}}_{Z}
\longrightarrow 
0,
$$
where the first map is defined by $f$, we have 
$H^i(Z, {\mathcal{O}}_Z)=0$  for all $i>0$   by the
corresponding long exact sequence and the  assumption
$H^i(Y, {\mathcal{O}}_Y)=0$  for all $i>0$. 

The regular function $f$ gives a rational map from $X$ to ${\Bbb{P}}^1$.
It is a morphism when restricted to $Y$. 
Let $C'=f (Y)$, the image of $Y$ under the map $f$.
 By Hironaka's elimination of indeterminacy,  
 there is a smooth  projective variety $X'$ such 
 that the morphism $\sigma : X'\rightarrow X$ is 
 composite of finitely many monoidal transformations.
$\sigma$  is an isomorphism  when restricted to $Y$.
So   $Y$ is fixed  and 
$g=f \circ \sigma: X'\rightarrow {\Bbb{P}}^1$ is proper,   surjective  and we have  
a commutative diagram
\[
  \begin{array}{ccc}
    Y                           &
     {\hookrightarrow} &
    X'                                 \\
    \Big\downarrow\vcenter{%
        \rlap{$\scriptstyle{g|_Y}$}}              &  &
    \Big\downarrow\vcenter{%
       \rlap{$\scriptstyle{g}$}}      \\
C'       & \hookrightarrow &
 {\Bbb{P}}^1.
\end{array}
\] 
Notice that the $D$-dimension does not change 
under blowing up or blowing down: 
Let  $\xi: V\rightarrow  W$
be a  surjective  morphism  of 
two varieties  and let  $D$ be  a Cartier divisor  on $W$,
then we have  (\cite{30}, Chapter 2,  Theorem 5.13)
$$  \kappa(\xi^*D,   V)=
\kappa(D, W). 
$$
In particular, in the above blowing up  $\sigma : X'\rightarrow X$, let 
$E$ be  an effective divisor on $X'$ such that
        codim$\sigma(E)\geq 2$, then  
        $$\kappa (\sigma^{-1}(D)+E, X')=\kappa (D, X), 
        $$
where  $\sigma^{-1}(D)$   is the reduced  transform of $D$,
defined to be
$\sigma^{-1}(D)=\sum D_i$, $D_i$'s are the   irreducible components
of  $D$.

The $D$-dimension  also does not depend on the choice of coefficients  
if   $D$  is an effective  divisor 
       with simple  normal crossings.
       Let $D_1$, $D_2$, $\cdot$$\cdot$$\cdot$, $D_n$
       be  any divisor on $X$ such that for every $i$, $1\leq i \leq n$, 
       $\kappa (D_i, X)\geq 0$, then for   integers 
       $p_1> 0,\cdot\cdot\cdot$, $p_n>0$,  
       we have   \cite{10}
       $$\kappa(D_1+\cdot\cdot\cdot+D_n,X)=
       \kappa(p_1D_1+\cdot\cdot\cdot+p_nD_n,X).
       $$
      In particular, if $D_i$'s are irreducible components of
      $D$ and $D$ is effective, then we can change its coefficients
      to different positive integers but  do not change the $D$-dimension.

A fiber space  is defined to be a  proper surjective   morphism
$f: V\rightarrow  W$  between  two varieties  $V$ and $W$
such that  the  general   fiber  is  connected. 
We cannot use the  above     morphism 
in the commutative diagram from
$X'$  to ${\Bbb{P}}^1$
 to  compute 
$D$-dimension because we do not know whether  a general  fiber is connected.

By Stein factorization, 
we can factor  the map $g$  through 
$$X'{\stackrel{h}{\longrightarrow}}\bar{C}
{\stackrel{\pi}{\longrightarrow}}
{\Bbb{P}}^1
$$
where   $g=\pi\circ h$,
$\pi$  is a  finite  morphism, 
 $h$  is proper, surjective  and every fiber   of  $h$  in  $X'$ is connected.
Let   $C=h(Y)$,  then we have the following new 
commutative diagram
\[
  \begin{array}{ccc}
    Y                           &
     {\hookrightarrow} &
    X'                                 \\
    \Big\downarrow\vcenter{%
        \rlap{$\scriptstyle{h|_Y}$}}              &  &
    \Big\downarrow\vcenter{%
       \rlap{$\scriptstyle{h}$}}      \\
C       & \hookrightarrow &
 \bar{C}.
\end{array}
\] 
Now we have a fiber space such that every fiber in
$X'$  is connected. 
Consider the image of $D$ under $h$. If $h(D)$   has dimension 0, 
then $Y$ contains  complete   curves, so $h(D)=\bar{C}.$

Since  $Z=\{f=0, f\in H^0(Y, {\mathcal{O}}_Y)\}$ 
is irreducible   and  $g|_Y=f$,  $\bar{Z}$
  is an irreducible component 
of  $g^{-1}(0)$ 
and  $g^{-1}(0)\cap  Y =\bar{Z}\cap  Y= Z$
is irreducible.

By a theorem of  Iitaka ( \cite{30}, Chapter II, Theorem 5.11), for 
a general point  $t$ on $C$,
we have 
$$ d= \kappa(D,   X')\leq  \kappa(D_t,   X'_t) + 1 \leq d,
$$
where   $X'_t=h^{-1}(t)$ is the 
irreducible 
fiber in $X'$
and $D_t$ is the restriction of  $D$  on $X'_t$.  
Thus $\kappa(D_t,   X'_t)=d-1$. 
Let  $t_0\in C$  such that  $\pi (t_0)=0$
and $\bar{Z}$  is an irreducible 
component of  $h^{-1}(t_0)=X_0$.
Then  $h^{-1}(t_0)$  is a connected 
component of  $g^{-1}(0)$.
By upper semi-continuity  theorem, 
$\kappa(D_0, X'_0)=d-1$,
where $D_0=D|_{X'_0}$. 
By the properties of $D$-dimension, 
$\kappa(D|_{\bar{Z}}, \bar{Z})=\kappa(D_0, X'_0)=d-1$.

\begin{flushright}
 Q.E.D. 
\end{flushright}

\begin{lemma} Under the condition of   Theorem  1.3, 
if  $Y$  is smooth, then 
  the regular functions on $Y$  separate points on $Y$. 
\end{lemma}
 $Proof.$ For two distinct points $y_1$ and $y_2$  on  $Y$,
we need to find a regular function $R$ on $Y$  such that
$R(y_1)\neq R(y_2)$.  We will use  induction on the dimension of  $Y$.
When $Y$  is a curve, the claim is true by Lemma  2.3. 
We may assume that the claim is true for $(d-1)$-dimensional 
varieties. Consider the  $d$-dimensional variety  $Y$.  

By Lemma 2.2, there is a smooth prime 
principle  divisor  $Z$  passing through 
$y_1$. Let $f\in  H^0(Y, {\mathcal{O}}_Y)$  be the defining  function of  
$Z$. 
If $f(y_2)\neq 0$, we are done. Assume $f(y_2)=0$. 
Consider the short exact sequence

$$ 0\longrightarrow 
 {\mathcal{O}}_Y
\longrightarrow 
 {\mathcal{O}}_Y
\longrightarrow 
 {\mathcal{O}}_{Z}
\longrightarrow 
0,
$$
where the first map is defined by $f$,  we have a surjective map
from $H^0(Y, {\mathcal{O}}_Y)$  to $H^0(Z, {\mathcal{O}}_{Z})$
since  $H^i(Y, {\mathcal{O}}_Y)=0$  for all  $i>0$.

By  Lemma  2.3  and  the inductive  assumption, 
$Z$  is affine. Therefore    
there is a regular function  $r$  on $Z$  such that
$r(y_1)\neq r(y_2)$. Lift this function from $Z$ to $Y$,
we get a global regular function $R$ on $Y$ 
such that it separates $y_1$  and $y_2$.

\begin{flushright}
 Q.E.D. 
\end{flushright}

Lemma 2.4 only holds for smooth varieties. 
The reason is that if $y_0$  is a singular point on
$Y$, we cannot expect that there 
exists a prime principal divisor passing through 
$y_0$. For example, let $A={\Bbb{C}}[x, y, z]/(xy-z^2)$
and $Y=$Spec$A$. However, to prove that
$Y$  is an  affine  variety, the following lemma 
is sufficient.

\begin{lemma} Under the condition of  Theorem 1.3, 
if  $Y$  is not smooth  but normal, then 
for any irreducible curve 
$F$  on  $Y$, there is a regular function 
$f$ on $Y$  such that the restriction
$f|_Y$  is  not  a constant.
\end{lemma}

$Proof$. Let  $\bar{F}$  be an  irreducible 
complete curve on  $X$ such that
$F$  is an open subset of  $\bar{F}$. 
Then  the  complement  $F^c=\bar{F}-F$
 is a set of finitely many  points 
 on  the boundary  $X-Y$. An old theorem of  Seidenberg 
 says that  a general hyperplane 
 section  $H$  of $X$  is normal and irreducible \cite{25}.
 We may choose $H$  such that all  points
 of $F^c$  are  not contained in  $H$.  Since  $H$  is ample and 
 $\bar{F}$ is complete, $H\cap \bar{F}=H\cap  F\neq  \emptyset$. 
 If  $F$  is contained in  $H$, then 
 $\bar{F}$  is a curve on  $H$  since  $\bar{F}$ is connected. So 
 $F$  is not contained in  $H$, i.e., 
 there is a point $p$  in  $F$  such that 
 $p$  is not a point of  $H$.

 Let  $h$
be the defining  linear form of 
$H$. Let $h'$ be a different linear 
form (linearly independent with  $h$)
such that the hyperplane $H'$ defined by 
$h'$  does not contain the point $p$. 
Then $h/h'$  gives a rational function 
on  $X$  and defines $H\cap X$.

Since  $D$  is a big divisor,
by  \cite{19},
there are two regular functions  $f$  and
$g$  on  $Y$  such that 
$h/h'=f/g$.   So  the irreducible open 
subvariety $Z=H\cap  Y$  is defined by  
$f$. Since  
$F$  is  not a curve on 
$Z$, the restriction function
$f|_F$  is not a constant. 

 \begin{flushright}
 Q.E.D. 
\end{flushright}

\begin{lemma} $Y$  is a quasi-affine  
variety
 under the assumption of  Theorem 1.3. 
\end{lemma}
$Proof$. By \cite{4},  there is a proper morphism  $\xi: Y\rightarrow U $
to a quasi-affine  variety  $U$
since 
for any irreducible curve 
$F$  on  $Y$, there is a regular function 
$f$ on $Y$  such that the restriction
$f|_Y$  is  not  a constant.
We know that  $Y$  has no complete curves, so
the fiber of the map $\xi$ is of  0 dimensional and
finite. Therefore  $\xi$  is a quasi-finite  morphism.
Zariski's Main Theorem (\cite{20},  Chapter  III, 
Section  9)  says that 
if  $\xi: Y\rightarrow U$ is  a 
morphism of varieties  with finite fibers,
then the map $\xi$  can be factored 
through an embedding $i$  from  
$Y$  to a variety  $Y'$  
followed by a finite morphism
$\psi: Y'\rightarrow  U$. 
Hence
$Y$  is a quasi-affine  variety
since $Y'$  is  quasi-affine. 

\begin{flushright}
 Q.E.D. 
\end{flushright}

{\bf Proof of Theorem  1.4.} 
In 1988, Neeman proved a very nice local
criterion for affineness: Let  $V=$${\mbox{Spec}}$A be a scheme  and 
$U\subset V$  a quasi-compact  Zariski  open  subset.
Here  we don't  assume that  $A$  is   noetherian 
\cite{23}. 
Then 
$U$  is affine  if and only if  $H^i(U, {\mathcal{O}}_U)=0$
for  all  $i\geq 1$.

Combining with 
 Theorem 1.3,  we immediately see that 
 $Y$
is affine.

\begin{flushright}
 Q.E.D. 
\end{flushright}

{\bf Proof of Theorem 1.5.} Let $X$  be an irreducible 
projective variety containing $Y$. Let   $V$ be 
an irreducible 
projective variety containing $W$. 
Then  we have a rational map $g$ from  $X$
to  $V$ such that  $g|_Y=f$. 
By Hironaka's elimination of indeterminancy,
we may assume that $g$ is a proper surjective morphism. 
Let   $\pi_V: V'\rightarrow V$
be  the blow up  of  closed subset of 
$V$ 
such that $V'$  is smooth.
Let  $\pi_X: X' \rightarrow  X $  be the 
resolution  of   the singularities of 
$X$  
such that   we have the following commutative diagram

\[
  \begin{array}{ccccc}
    Y                          &
     {\hookrightarrow} &
    X     &   {\leftarrow} & X'                       \\
    \Big\downarrow\vcenter{%
        \rlap{$\scriptstyle{g|_Y=f}$}}         &&  
    \Big\downarrow\vcenter{%
       \rlap{$\scriptstyle{g}$}}    &&
    \Big\downarrow\vcenter{%
       \rlap{$\scriptstyle{h}$}}    \\
W       & \hookrightarrow &
V &\leftarrow&  V' ,
\end{array}
\]
where $f$ is a surjective morphism, $g$  and $h$  are   
  surjective proper morphisms.

Suppose that the dimension of  
$Y$  is $d$  and the dimension of  
$W$  is  $m$.   
Choose suitable 
projective variety   $V$
such  that the boundary $V-W$
is support of an ample divisor
  $A$.
  Let  $D$  be an effective 
divisor on $X$  with  support  $X-Y$.
Then we have  (\cite{30}, Chapter II, Theorem 5.13)
$$\kappa(\pi_V^*A, V')=\kappa(A, V)=m$$
and  $$\kappa(\pi_X^*D, X')=\kappa(D, X).$$ 

If we can prove that  $Y$  has no complete curves and 
$\kappa(D, X)=d$, then  $Y$  is affine by Theorem 1.4. 
The first property is obvious since every fiber  of  
$f$ 
in  $Y$  is affine. To compute the $D$-dimension of 
$X$, we need a theorem of Fujita:  
 Let 
$M$ and  $S$  be two projective manifolds.  
 Let $\pi: M\rightarrow S$
  be a  fiber space and let $L$ and $H$ be  line bundles  on $M$
  and $S$ respectively. Suppose that  
  $\kappa (H, S)=\dim S$ and that
  $\kappa (aL-b\pi^*(H))\geq 0$
for certain  positive integers $a$, $b$.  Then
$$\kappa(L,M)=\kappa (L|_F, F)+ \kappa (H, S)$$
for a general  fiber $F$  of  $\pi$.   
Here if $L$  is a  line  bundle  
on a  projective  manifold   $M$, 
it determines   a  Cartier  divisor  $D$.
We  define   $$\kappa(L, M)=\kappa(D,  M).$$

 We don't know whether the fiber  of 
$h$ is connected. Let $$X'{\stackrel{j}{\longrightarrow}}V''
{\stackrel{\alpha}{\longrightarrow}}
V'
$$
be the Stein factorization, then 
$h=\alpha \circ j$, $\alpha$  is a finite morphism
and   $j$  has connected fibers. And a general fiber
of $j$  is smooth and irreducible.

 Since  the support of  $D$ is $X-Y$, the image 
 $\pi_V\circ h(\pi_X^*D)=V-W$ is support of 
 $A$.  So  the support of pull back divisor $h^*(\pi^*_V(A))$  is  contained in 
 the support of  $\pi_X^*(D)$.

 Since  the dimension of  $V''$  is  $m$, we have 
 (\cite{30}, Chapter II, Theorem 5.13)
 $$\kappa(\alpha^*(\pi_V^*(A)), V'' )=\kappa(A, V)=m.
 $$
Let  $H=\alpha^*(\pi_V^*(A))$, then $\kappa(H, V'')=m=$dim$V''$. Let  $L=\pi_X^*D$, then 
for sufficiently large $n$, we have 
$$  \kappa(nL-j^*H)\geq 0.
$$
This is because the exceptional divisors do not
change the  Iitaka dimension and 
the support of  $g^*(A)$  is contained in  $X-Y$,
the support of  $D$. 
By Fujita's formula, 
$$\kappa(L, X')=\kappa(nL, X')=\kappa(nL|_F, F)+m,
$$
where $F$  is a general fiber of  $j$. 
A general fiber of $g$  has dimension 
$d-m$. 
Since  every fiber of  $f$  is affine and $\pi_X(F)$
is a fiber  of  $g$, we have 
$$\kappa(nL|_F, F)=\kappa(nD|_{\pi_X(F)}, \pi_X(F))=
d-m. $$
Therefore 
$$\kappa(D, X)=\kappa(L, X')=\kappa(nL|_F, F)+m=d.
$$
Hence  $Y$  is affine by Theorem 1.4. 
\begin{flushright}
 Q.E.D. 
\end{flushright}

{\bf Proof of Corollary 1.6.}
By the same calculation as in the proof of 
Theorem 1.5, we have  $\kappa(D, X)=d$, the dimension of  $Y$. 
The conclusion is obvious by Theorem 1.4.

\begin{flushright}
 Q.E.D. 
\end{flushright}

\section{Examples}

Again $Y$  is a variety  contained in a projective variety
$X$  such that  $Y=X-D$, where $D$  is an effective  boundary 
divisor   with  support  $X-Y$.

\begin{example} There is  an  affine surface $Y$ such that 
the  graded ring 
$\oplus_{n= 0}^\infty  H^0(X, {\mathcal{O}}_X(nD))$  is not finitely 
generated.  This example is due to  Zariski (\cite{31}, page 562-564).

Let  $C$  be a  smooth  curve of
degree 3  in  ${\Bbb{P}}^2$.  Let  
$\Lambda$ be a divisor  class cut  out  on  $C$
by  a curve  of  degree  4  in 
${\Bbb{P}}^2$.  There exist  12 
distinct  points  $p_1$, $p_2$, $\cdot\cdot\cdot$, $p_{12}$
on $C$  such that
$$  m(p_1+p_2+\cdot\cdot\cdot+p_{12})\not\in m \Lambda
$$
for  all  positive integers  $m$.  
Let  $X$  be  a surface  obtained by  blowing  up 
${\Bbb{P}}^2$  at these 12 points   
$p_1$, $p_2$, $\cdot\cdot\cdot$, $p_{12}$. 
Let  $\bar{C}$  be the strict  transform  of  $C$
(i.e., the closure  of the inverse image  of  
$C-\{p_1,p_2,\cdot\cdot\cdot,p_{12}\}$  in  $X$).
Let  $L$
be a line  not passing  through  any  point  $p_i$
in these 12 points.  Let  $\bar{L}$  be the strict 
transform of  $L$. 
Then  the complete  linear system 
$$|m(\bar{C}+\bar{L})|$$
    has a  fixed locus  
$\bar{C}$  for all  $m\geq 1$  and  
$$|m\bar{C}+(m-1)\bar{L}|$$
has no   fixed  components and  is  base point free. 
By Nakai-Moishenzon's  ampleness 
criterion  (\cite{7}, Chapter  V, Section 1),
the divisor 
$$m\bar{C}+(m-1)\bar{L}$$
   is ample.
Hence  the complement  $Y=X-(m\bar{C}+(m-1)\bar{L})$
is affine  but the  graded ring 
$$R=\oplus_{m= 0}^\infty H^0(X, {\mathcal{O}}_X(m(\bar{C}+\bar{L})))$$
is  not finitely  generated.

\end{example}

\begin{example}  There is a   nonaffine surface $Y$  such that 
the  graded ring 
$$\oplus_{n= 0}^\infty H^0(X, {\mathcal{O}}_X(nD))$$
  is  finitely 
generated  \cite{16}. 

Let $C$ be an elliptic curve and $E$ the unique nonsplit 
      extension of $\mathcal{O}$$_C$ by itself.  
      Let ${X=\Bbb{P}}_C(E)$ and  $D$ be the canonical section, then 
$Y=X-D$  is not affine    and   $H^0(X, {\mathcal{O}}_X(nD))=\Bbb{C}$ \cite{16}.
So  $$\oplus_{n= 0}^\infty H^0(X, {\mathcal{O}}_X(nD))$$  is  finitely 
generated. 
\end{example}
 The above two examples demonstrate that the  affineness of $Y$
and the finitely  generated  property of  the  graded ring 
$$\oplus_{n= 0}^\infty H^0(X, {\mathcal{O}}_X(nD))$$  
are different in nature. 
The reason is  
$$\Gamma(Y, {\mathcal{O}}_Y )\neq 
\oplus_{n= 0}^\infty H^0(X, {\mathcal{O}}_X(nD)).$$ 
In fact,  we have  \cite{4}

\begin{lemma}{\bf[Goodman, Hartshorne]}  Let $V$ be a scheme and  
$D$ be an effective Cartier divisor on  $V$. Let $U=V-$Supp$D$ and $F$ 
be any coherent sheaf on $V$, 
then for every $i\geq 0,$ 
$$\lim_{{\stackrel{\to}{n}}}
H^i(V, F\otimes {\mathcal{O}}(nD)) \cong  H^i(U,  F|_U).
$$
\end{lemma}

So  we have 
$$ \Gamma(Y, {\mathcal{O}}_Y )\cong
\lim_{{\stackrel{\to}{n}}}
H^0(X,  {\mathcal{O}}_X(nD)).
$$

The direct limit is  the quotient of 
the direct sum and its subring, so it is
 much ``smaller"  than direct sum (\cite{17}, Chapter II, Section  10). 
And  even though  $Y$  is affine, the boundary divisor can be very 
bad. For example, $D$ may not be nef. It is easy to see this by blowing up 
${\Bbb{P}}^2$ at a point. Let $L$ be a line in 
${\Bbb{P}}^2$, let $O$  be a point on $L$. Let 
$\pi: X\rightarrow {\Bbb{P}}^2$  be the
blow up  of 
${\Bbb{P}}^2$ at  $O$. Let $E$  be the exceptional  divisor
and  $D=\pi^{-1}(L)+ mE$, where  $\pi^{-1}(L)$
is the strict transform of  $L$  and $m$  is a large 
positive integer.  Then  $D\cdot E=1-m <0$  (\cite{7}, 
Chapter V, Corollary 3.7).  Therefore $D$ is not nef.

\begin{example} There exists  a   threefold  $Y$ such that  $Y$  contains no complete
curves, $H^i(Y, \Omega^j_Y)=0$  for all $i>0$  and $j\geq 0$  but is not affine
\cite{34}.

Let $C_t$ be a smooth projective elliptic curve 
defined by $y^2=x(x-1)(x-t)$, $t\neq 0, 1$.  Let $Z$ be the elliptic surface
defined by the same equation, then we have surjective morphism from $Z$ to 
$C={\Bbb{C}}-\{0, 1\}$ such that for every $t\in C$, the fiber $f^{-1}(t)=C_t$.   
In \cite{34}, we proved that
there is a rank 2 vector bundle $E$ on $Z$ such that when 
   restricted to $C_t$, $E|_{C_t}=E_t$ is the unique nonsplit 
   extension of 
    ${\mathcal{O}}_{C_t}$ by ${\mathcal{O}}_{C_t}$, where $f$
     is the morphism from $Z$ to $C$. 
 We also proved that 
   there is a divisor $D$ on 
$X={\Bbb{P}}_Z(E)$ such that when restricted to $X_t={\Bbb{P}}_{C_t}(E_t)$,
$D|_{X_t}=D_t$ is the canonical section of $X_t$.
Let $Y=X-D$,
 we have $H^i(Y, \Omega^j_Y)=0$
for all $i>0$ and $j\geq 0$  \cite{31}.  
We know that  this  threefold  $Y$  contains no complete curves \cite{33}
and  $\kappa(D, X)=1$.  So $Y$  is not affine.   
\end{example} 

\begin{example} There is a  surface $Y$  without complete curves  such that 
$\kappa(D, X)=2$  but  is not affine.

 Remove a  line  $L$  from  ${\Bbb{P}}^2$, we have 
${\Bbb{C}}^2={\Bbb{P}}^2-L$. Remove the origin  $O$  from 
the complex plane  ${\Bbb{C}}^2$, let  $Y={\Bbb{C}}^2-\{O\}$.
Then  $Y$  is not affine since the boundary 
is not connected  (\cite{8}, Chapter II, Section  3 and Section  6).
Blow up   ${\Bbb{P}}^2$  with center $O$, let $E$
be the exceptional  divisor  and $\pi: X\rightarrow {\Bbb{P}}^2$
be the blowup. Let $D=\pi^{-1}(L)+E$, where  $\pi^{-1}(L)$
is the strict transform of  $L$. Then 
by Iitaka's  result, on $X$, $\kappa(D, X)=2$ and  $X-D\cong  Y$
has no complete curves. But $Y$  is not affine. 
\end{example}

\begin{lemma}{\bf[Koll$\mbox{\'{a}}$r]} 
Let  $\pi: X\rightarrow Z$
be  a  surjective  map between  projective   varieties,
$X$ smooth, $Z$  normal.  Let
$F$  be the  geometric  generic  fiber  of
$\pi$  and  
assume that  
$F$ is connected.  
The following two statements   are equivalent:

(i)  $R^i\pi_*{\mathcal{O}}_X=0$  for all  $i>0$;

(ii)  $Z$  has  rational  singularities  and  
$H^i(F, {\mathcal{O}}_F)=0$
for all $i>0$. 
\end{lemma}

\begin{example} There is a smooth  variety $Y$  of dimension  $d\geq 1$ 
with  $H^i(Y, {\mathcal{O}}_Y)=0$
for all $i>0$  and  $\kappa(D, X)=d$  but  is not affine.

Let  $X$  be the smooth projective  variety 
    obtained by  blowing up   a  point  $O$
in  ${\Bbb{P}}^d$. Let  $\pi: X\rightarrow  {\Bbb{P}}^d$ be the blowup. 
Let  $H$  be a hyperplane  not passing through $O$. 
 Let $D=\pi^{-1}(H)$, the strict transform of 
 $H$    and
$Y=X-D$.  Then   $\kappa(D, X)=d$  (\cite{30}, Chapter 2, Theorem 5.13).  
Let $E$ be the exceptional divisor on  $X$, then 
$E\cong {\Bbb{P}}^{d-1}$. So  
$H^i(E, {\mathcal{O}}_E)=0$
for all $i>0$. 
By   Lemma 3.6, 
we  have   $R^i\pi_*{\mathcal{O}}_X=0$  for all  $i>0$.

We can get the vanishing of direct images 
by Grauert's upper semi-continuity 
theorem. 
For every point $p\in {\Bbb{P}}^d$, $p\neq O$, then
$\pi^{-1}(p)=q$  is a point in  $X$. So
for every $i>0$, 
$$ h^i(q, {\mathcal{O}}_q)=0.
$$
We already saw 
$$h^i(\pi^{-1}(O), {\mathcal{O}}_{\pi^{-1}(O)})=0.
$$
So  $R^i\pi_*{\mathcal{O}}_X=0$  for all  $i>0$
(\cite{30}, Chapter 1, Theorem 1.4).

Let $U={\Bbb{P}}^d-H$, then $U\cong {\Bbb{C}}^d$.
 For the global sections on affine  space $U$,  we have
(\cite{7}, Chapter III,   Proposition 8.1, 8.5 and Chapter II, Proposition 5.1(d))
$$0=R^i\pi_*{\mathcal{O}}_X(U)= H^i(Y, {\mathcal{O}}_Y) 
$$
for all $i>0$.

It is obvious that $Y$  is not affine since it contains a projective space
${\Bbb{P}}^{d-1}$.

\end{example}
\begin{remark}
Examples 3.5-3.7  show  that  $Y$  is not affine
if we drop any condition  from three conditions in  Theorem 1.4.
\end{remark}

\begin{example} There is a  threefold  $Y$  such that it satisfies the followings three conditions but is not affine:

(1) $Y$  contains no complete curves;

(2) the boundary  $X-Y$  is  connected;

(3) $\kappa(D, X)=3$.

Let $H$  be a hyperplane in  ${\Bbb{P}}^3$.
Let $L$  be a line not contained in $H$. 
Blow up ${\Bbb{P}}^3$ along $L$, let
$\pi: X\rightarrow  {\Bbb{P}}^3$  be the blowup.
Define a divisor $D$ on $X$ such that 
$D=\pi^{-1}(H)+E $, where $E$  is the exceptional 
divisor. Let  $Y=X-D$, then $Y={\Bbb{P}}^3-H-L. $ 

It is easy to see that the above three 
conditions are satisfied but $Y$  is not affine.  
\end{example}

\begin{remark}

Theorem 1.2 is false for 
higher dimensional varieties by the above example. 
\end{remark}

\begin{example} Let $(x, y)$ be the coordinates of 
${\Bbb{C}}^2$. Define a projection 
$$\pi: {\Bbb{C}}^2-\{(0,0)\}\longrightarrow \Bbb{C} $$
$$\quad\quad (x,y)\longmapsto x.$$
Then $\pi$ is a surjecive morphism  and  every fiber is affine. 
Let
$Y={\Bbb{C}}^2-\{(0,0)\}$, then the boundary 
of  $Y$  in  ${\Bbb{P}}^2$
is not connected. Therefore
 ${\Bbb{C}}^2-\{0\}$  is not affine (\cite{8}, Page 67). 
By Neeman's theorem \cite{23}, 
we know $H^i(Y, {\mathcal{O}}_Y) \neq  0$.  
Thus if we drop the cohomology 
assumption, then Theorem 1.5 does not hold. 
\end{example}

\begin{example}  There is a quasi-projective variety $Y$  with
a surjective morphism  $f: Y\rightarrow U$  such that
$U$  is affine, a general fiber is affine and 
$H^i(Y, {\mathcal{O}}_Y)=0$ for all $i>0$  but $Y$  is not affine.

Let  $Y$, $U$  be the varieties defined in 
Example 3.7. Then the fiber space $\pi: Y\rightarrow U$ 
satisfied the above requirements. $Y$  is not affine 
because it has a projective space 
${\Bbb{P}}^{d-1}$.

\end{example}

\end{document}